\documentclass{article}

\usepackage{amsmath, amssymb, graphicx, amsthm}

\def\ds{\displaystyle}

\def\veps{\varepsilon}

\def\R{\mathbb{R}}
\def\Z{\mathbb{Z}}
\def\N{\mathbb{N}}

\def\Z{\mathbb{Z}}

\def\sign{{\rm{sign\,}}}

\newtheorem{theorem}{Theorem}[section]
\newtheorem{example}{Example}[section]
\newtheorem{proposition}[theorem]{Proposition}
\newtheorem{remark}[theorem]{Remark}
\newtheorem{lemma}[theorem]{Lemma}
\newtheorem{definition}[theorem]{Definition}

\title{Global Solvability in Functional Spaces for \\ Smooth Nonsingular Vector Fields in the Plane}

\author{Roberto De Leo $\langle$deleo@unica.it$\rangle$, Todor Gramchev $\langle$todor@unica.it$\rangle$\\
    Dipartimento di Matematica e Informatica\\ 
    Universit\`a di Cagliari, \\
    Via Ospedale 72 - 09124 Cagliari, Italy
    \and 
    Alexandre Kirilov $\langle$akirilov@ufpr.br$\rangle$\\
    Departamento de Matem\'atica \\
    Universidade Federal do Paran\'a \\
    Caixa Postal 19081\\
    81531-990 Curitiba, Paran\'a, Brasil
    \thanks{The third author was partially supported by CNPq, Brasil.}
}
\begin{document}
\maketitle
\begin{abstract}
We address some global solvability issues for classes of smooth nonsingular vector fields $L$ in the plane related to cohomological equations $Lu=f$ in geometry and dynamical systems. The first main result is that $L$ is not surjective in $C^\infty(\R^2)$ iff the geometrical condition -- the existence of separatrix strips -- holds. Next, for nonsurjective vector fields, we demonstrate that if the RHS $f$ has at most infra-exponential growth in the separatrix strips we can find a global weak solution $L^1_{loc}$ near the boundaries of the separatrix strips. Finally we investigate the global solvability for perturbations with zero order p.d.o. We provide examples showing that our estimates are sharp.


\end{abstract}

%%% ----------------------------------------------------------------------
%%% ----------------------------------------------------------------------
%\tableofcontents

\section{Introduction and main results} \

We recall that Duistermaat and H\"ormander, see \cite{DH1}, have demonstrated that a nonsingular smooth vector field $X$ in a $n$-dimensional open manifold $M$ is surjective if and only if it admits a global transversal section, namely a smooth hypersurface which is transversal to $X$ at every point and cuts exactly once every of its integral trajectories.

On the other hand, nonsurjective vector fields appear in the context of the geometry of foliations (see \cite{Re1}) and dynamical systems ( see \cite{Del1,Nov1}). In particular, the issue of the global solvability of the cohomological equations of the type $Xu =f$ is a challenging and difficult problem related to Geometry, Dynamical Systems( cf. \cite{Fo1}, see also \cite{Grom1} on the solvability of systems of PDEs) and in the general theory of PDEs, e.g. see \cite{BSi,BKi,GDY,PeZa1} on global solvability on tori, and \cite{CGR1,DW1,GPR1} in the Gelfand--Shilov spaces $S^\mu(\R^n)$. Finally, we mention that the surjectivity in various functional spaces for linear partial differential operators of higher order have been extensively studied since 80's (see \cite{AlPo1,BMT1,BMT2} and the references therein).

In this work we investigate, in the framework of the general theory of PDEs, the global solvability in the plane for smooth nonzero vector fields which appear in theory of foliations and the cohomological equations in Geometry and Dynamical Systems. We also investigate the stability of the global solvability in weighted Sobolev spaces under perturbation with zero order pseudodifferential operators.

We consider smooth nonsingular real vector field in he plane
\begin{equation}
Lu  =  p(t) \partial_t u + q(t)\partial_xu = f(t,x),
\label{vf3}
\end{equation}
i.e.,  $p$ and $q$ are real-valued smooth functions which have no common zeros.
%\begin{equation}
%|p(t)| + |q(t)| \neq 0, \quad t\in \R.
%\label{vf3nz}
%\end{equation}

One assumes that there is an integer $N\geq 2$ and $t_1 < \ldots < t_N$ such that
\begin{equation}
p(t) =0 \Longleftrightarrow t=t_j, \, j=1,2, \ldots , N
\label{vf4}
\end{equation}
with
\begin{equation}
p'(t_j) \neq 0, \quad  \, j=1,2, \ldots , N
\label{vf4a}
\end{equation}
and
\begin{equation}
\textrm{ $q$ admits at most one zero in $]t_j, t_{j+1}[$ for $j=1,2, \ldots , N-1$.}
\label{vf4b}
\end{equation}

\smallskip
Note that the lines $\{t=t_j\}$, $j=1,\ldots, N$, are characteristics for $L$. We also suppose that $p$ and $q$ are polynomials. 

Our results are true under weaker restrictions on $p$ and $q$, but we prefer to exhibit the main novelties avoiding highly technical arguments and capturing particular cases of $L$ of interest in geometry and dynamical systems
(cf. \cite{Ca1} for foliations, see also \cite{Del1} for a thorough discussion of its action on $C^\infty(\R^2)$). For example,
\begin{equation}
L_0 u  =  (1-t^2) \partial_t u - 2t \partial_xu
\label{vf5}
\end{equation}
and more generally,
\begin{equation}
L_{\lambda,k}u  =  (1-t^2) \partial_t u + \lambda t^{k} \partial_xu
\label{vf5a}
\end{equation}
for $\lambda \neq 0$, $k\in \N$.

\smallskip
The first main goal of the present work is to show that the existence  of separatrix type phenomena for \eqref{vf3} is the only obstruction for the surjectivety in $C^\infty(\R^2)$ of $L$. Moreover, we exhibit functional spaces associated to the separatrix strips where we can solve globally this cohomological equation in $\R^2$ and investigate the stability of this global solvability under perturbations of $L$ with  zero order pseudodifferential operators in $x$.

\begin{definition}
A strip
\begin{equation}
S_j  =  \{(t,x): \, t\in ]t_j, t_{j+1}[,  x\in \R  \}, \ \mbox{ with } j\in \{1, \ldots, N-1 \}
\label{sep1}
\end{equation}
is a separatrix for the vector field $L$ above if all characteristic curves $x=x(t; \tau, y)$, starting at a point $(\tau, y) \in S_j$ satisfy
\begin{align}
 \textrm{either }\ \lim_{t\to t_j^+} x(t;\tau, y) = \lim_{t\to t_{j+1}^-} x(t;\tau,y) = +\infty,\\
 \textrm{or } \ \lim_{t\to t_j^+} x(t;\tau, y)= \lim_{t\to t_{j+1}^-} x(t;\tau, y) = -\infty.
\label{sep2}
\end{align}
\end{definition}

We state the first new result of our article.

\begin{theorem}\label{Theorem1.1}
The following assertions are equivalent:
\begin{enumerate}
\item [ i)] the vector field $L$ is not surjective in $C^\infty (\R^2)$;
\item [ ii)] the vector field $L$ admits a separatrix $S_j,$ for some $j\in \{1,\ldots, N-1 \}$;
\item [ iii)] there exists $j\in \{1,\ldots, N-1 \}$ and $\theta_j \in ]t_j, t_{j+1}[$ such that
$q(\theta_j) = 0$ and \\ $q$ has opposite signs in $]t_j, \theta_j[$ and $]\theta_j, t_{j+1}[$.
\end{enumerate}
In particular, the operators $L_{\lambda, k}$ are not surjective in $C^\infty (\R^2)$ if and only if $k$ is odd.
\end{theorem}

To illustrate the nonsurjectivity for simple example we point out that nonzero constants do not belong to $L_0(C^\infty (\R^2))$. Direct calculations implies that $L_0u=~c$ has a weak solution $u(t,x) =  \frac{c}{2} \ln \frac{|1+t|}{|1-t|}$. We show for more general classes of RHS $f\in C^\infty (\R^2)$ that  every solution has singularity either at $t=1$ or $t=-1$ (see Section 4 for more details).

This example shows that in order to solve globally $Lu = f$ one should allow some (weak) singularities of the type $L^1_{loc}$ near the adjacent characteristics forming the separatrix strips.

The second main novelty of this work is that, in order to find a global weak solution, in general the RHS $f(t,x)$ should grow at most like $\ds O (e^{\veps |x|}),$ for $|x| \to \infty$ uniformly in the separatrix strips $S_j$.

Finally, we derive sharp estimates on the singularities of the global solutions $u(t,x)$ of \eqref{vf3} near $t_j$, $j\in I_L$ for large classes of smooth RHS $f$, where \smallskip
\begin{equation}
I_L = \{ t_j:\, \textrm{$S_j$ or $S_{j-1}$ is  separatrix}, j=1, \ldots, N \}.
\label{sep2a}
\end{equation}

\smallskip
We point out that the part ii) of Theorem \ref{Theorem1.1} implies that $L$ is not surjective in $C^\infty (\R^2)$ if and only if $I_L$ is not empty. \\

In order to state the main result on the global solvability of \eqref{vf3} we introduce the subspace of the functions of infra-exponential growth in the $x$ variable (e.g., cf. \cite{Kan} where such growth plays an important role in theory of Fourier transform for hyperfunctions).
\begin{align}
C^\infty (\R: Exp_{sl} (\R)) & \doteq &\{ f\in C^\infty (\R^2): \, \forall T>0, \forall \veps >0, \forall \alpha \in \Z_+^2,   \exists C>0  \nonumber \\
& & \textit{ s.t. } |\partial_{t,x}^\alpha f(t,x) |\leq Ce^{\veps |x|}, \, |t|\leq T, x\in \R\} \label{subexp1}
\end{align}

We recall also the weighted Sobolev spaces $H^{s_1,s_2} (\R^n)$ in $\R^n$ (e.g. see \cite{Co1}).
\begin{equation}
H^{s_1,s_2}(\R^n)  \doteq \{ f\in \mathcal{S}'(\R^n):  \| f\|_{s_1,s_2} = \| \langle x \rangle ^{s_2} \langle D \rangle^{s_1} f \|_{L^2} < +\infty \}  \label{Sobc1}
\end{equation}
which measure the global regularity and the behaviour on $\infty$ in $\R^n$, where \linebreak $\langle x \rangle = \sqrt{1+\|x\|^2}$.\\

\begin{theorem}\label{Theorem2}
Let $L$ defined above be nonsurjective in $C^\infty(\R^2)$. Then we can find a right inverse $L^{-1}$ of $L$  acting continuously
\begin{equation}
L^{-1}: C^\infty (\R: Exp_{sl} (\R)) \longrightarrow   L^1_{loc} (\R: Exp_{sl} (\R))\, \mbox{$\bigcap$}\, C^\infty (\R\setminus I_L:Exp_{sl}(\R)) \label{inv1}
\end{equation}
and
\begin{equation}
L^{-1}: C(\R: H^{s_1,s_2} (\R)) \longrightarrow  L^1_{loc} (\R:  H^{s_1,s_2} (\R))\, \mbox{$\bigcap$}\,  C(\R\setminus I_L: H^{s_1,s_2}(\R)), \label{inv2}
\end{equation}
with $ s_1,s_2\in \R.$

Moreover, for any $\veps >0$ we have
      \begin{equation}\label{sup-general}
        \sup_{t \in [-\theta,\theta]}\left( \prod_{j=1}^{N}|t-t_j|^\veps \| L_j^{-1} f(t,\cdot) \|_{H^{s_1,s_2} (\R)}  \right) \leq C_{\veps, s_1, s_2, \theta} \| f \|_{ C(\, \overline{I}_j: H^{s_1,s_2} (\R))}
      \end{equation}

Next, if $f$ is a polynomial function with respect to $x$, i.e., $ f(t,x) = \sum_{\ell=0}^k f_\ell(t) x^\ell$, then
\begin{equation}
L^{-1} f(t,x) = \sum_{\ell=0}^k g_\ell(t) x^\ell \label{inv3}
\end{equation}
with
\begin{eqnarray}
  g_k(t) &=& O(\ln^{k+1} |t-t_j|) \textrm{ near $t= t_j$, if $S_j$ or $S_{j-1}$ is a separatrix} \label{inv4} \quad \\
  g_\ell(t) &=& o(\ln^{k+1} |t-t_j|) \textrm{ near $t= t_j$, if $S_j$ or $S_{j-1}$ is a separatrix,} \label{inv5}
\end{eqnarray}
for $\ell =0,\ldots, k-1.$ \\

Finally, given a zero order p.d.o. $b(t,x,D)$ in $x$ smoothly depending on $t$, and $s_1,s_2\in \R$ we can find $\veps_0 = \veps_0(L,s_1,s_2) >0$ such that if
\begin{equation}
\max_{{|\alpha|\leq [s_1]+2}\atop{|\beta| \leq [s_2]+2}}\sup_{{t\in [t_j,t_{j+1}]}\atop{(x,\xi)\in \R^2}} \langle x \rangle^{-\alpha} \langle \xi \rangle^{-\beta} |\partial_x^\alpha\partial_\xi^\beta b(t,x,\xi)|   < \veps_0
\end{equation}
then $L+b(t,x,D)$ admits a right inverse which satisfies \eqref{inv2}.\
\end{theorem}

The paper is organized as follows. Section 2 deals with the proof of Theorem~\ref{Theorem1.1} and exhibits some geometric features. We derive in the Section 3 precise estimates on suitable right inverses in the separatrix strips and proof a crucial gluing lemma. In Section 4 we obtain sharp results for $L_0$ on the singular behaviour near the separatrix lines. In Section 5 we consider perturbations of the nonsurjective vector field $L_0$ with a constant p.d.o. Finally, we discuss some possible generalizations in Section 6.

\section{Separatrix Strips and Nonsurjectivity}\

In this section we prove Theorem \ref{Theorem1.1}. We start by calculating the global ``singular" characteristics of $L$ after dividing by $p(t)$, namely, rewriting formally
 $Lu+bu = f$ to
\begin{eqnarray}
 \tilde{L}u +  \frac{1}{p(t)}b(t,x,D) u & = &  \frac{f(t,x)}{p(t)}
  \label{ops1}
\end{eqnarray}
with
\begin{eqnarray}
 \tilde{L}u  & = &  \partial_t u + \frac{q(t)}{p(t)} \partial_x u
  \label{ops2}
\end{eqnarray}

 The characteristics of $\tilde{L}$,  different from $t=t_j$, $j=1,\ldots, N$, are  defined by
 \begin{eqnarray}
 \dot{x} (t)  & = &  \frac{q(t)}{p(t)}, \qquad x|_{t=\tau} = y
  \label{ops3}
\end{eqnarray}
for some $\tau\neq t_j$, $j=1,\ldots, N$.\\

 We have

 \begin{lemma}\label{primitive-of-p-over-q}
The function $q(t)/p(t)$ has a global primitive $\rho(t)$ such that
\begin{eqnarray}
 \rho(t) & = &  \sum_{j=1}^N \varkappa_j q(t_j)\ln |t-t_j|  + \tilde{\rho}(t)
  \label{chars1}
\end{eqnarray}
where each $\varkappa_j\in \R\setminus \{ 0 \}$, with $j=1,\ldots, N$, depends only on $p(t)$ and $\tilde{\rho}\in C^\infty (\R)$.\\
Moreover, for each $j\in \{ 1, \ldots, N-1\}$ fixed, we have
\begin{eqnarray}
  \varkappa_j \varkappa_{j+1}q(t_j) q(t_{j+1}) >0  & \Leftrightarrow &  \textrm{$q$ admits a zero in $]t_j,t_{j+1}[$ of odd order} \qquad
  \label{chars2a}\\
  \varkappa_j \varkappa_{j+1}q(t_j) q(t_{j+1}) <0  & \Leftrightarrow &  \textrm{$q$ does not admit zero of odd order}
  \label{chars2b}
\end{eqnarray}

\end{lemma}

\begin{proof}
By the hypotheses \eqref{vf4}, \eqref{vf4a} on $p$ and the decomposition of rational functions, there are nonzero real numbers $\varkappa_1, \ldots, \varkappa_N$ and $r_1 \in C^\infty (\R)$ such that
\begin{eqnarray}
  \frac{1}{p(t)} & = &  \sum_{j=1}^N \frac{\varkappa_j}{t-t_j} + r_1(t)
  \label{chars3}
\end{eqnarray}
which yields
\begin{eqnarray}
  \frac{q(t)}{p(t)} & = &  \sum_{j=1}^N \frac{\varkappa_jq(t_j)}{t-t_j} + r_2(t)
  \label{chars4}
\end{eqnarray}
for some $r_2\in C^\infty (\R)$. The expression \eqref{chars1} follows by integration.

\smallskip
We note that the hypothesis \eqref{vf4} implies $q(t_j) \neq 0$, and hence
\begin{eqnarray}
  c_j & \doteq &  \varkappa_jq(t_j)\neq 0, \qquad j=1,\ldots, N
  \label{chars4a}
\end{eqnarray} \vskip-0.2em
\end{proof}

Next, we present an important auxiliary result.

\begin{lemma}\label{theta-equation}
Let $x(t, y)$ be defined by
\begin{equation}
  \dot{x} =  \frac{\lambda (t-\theta)^k}{(\theta_+ -t) (t-\theta_-)} + \tilde{q}(t), \quad x(\theta) = y, \quad \theta \in ]\theta_-, \theta_+[,
  \label{char1}
\end{equation}
with $\tilde{q}\in C^\infty ([\theta_-, \theta_+])$.

Then one can find $r\in C^\infty ([\theta_-, \theta_+])$ such that
\begin{equation}
x(t,y) = y + c_+ \ln |t-\theta_+| + c_- \ln |t-\theta_-| + r(t), \label{char2}
\end{equation}
where
\begin{equation}
c_\pm = \mp \frac{\lambda (\theta_\pm -\theta)^k}{\theta_+-\theta_-}  \label{char2-constants}
\end{equation} \\

\noindent In particular, we observe that
\begin{enumerate}
   \item[{\it i)}] $c_+c_- >0 \ \Leftrightarrow $ k is odd  $\Leftrightarrow \ c_+ $ and $c_-$ have the same signal and $\lambda > 0$;
   \item[{\it ii)}] $c_+c_- <0 \ \Leftrightarrow \ $ k is even $\Leftrightarrow \ c_+$ and $c_-$ have different signals and $\lambda < 0$.
\end{enumerate}
\end{lemma}

\begin{proof}
The proof follows from the decomposition
\begin{equation}
\frac{\lambda (t-\theta)^k}{(\theta_+ -t) (t-\theta_-)}
   =   \frac{\lambda (\theta_+ -\theta)^k}{(\theta_+-\theta_-) (\theta_+-t)}+
  \frac{\lambda (\theta_\pm -\theta)^k}{(\theta_+-\theta_-)(t-\theta_+)} + \tilde{q}_1(t),
  \label{char2a}
\end{equation}
where $\tilde{q}_1 = 0$ if $k=0,1$, and $\tilde{q}_1$ is polynomial of degree $k-2$, if $k\geq 2$,
and integration (from $\theta$ to $t$) of the RHS of \eqref{char1}.
\end{proof}

Now we present the main steps of the proof of Theorem \ref{Theorem1.1}. First, assume that $S_j$ is a separatrix, for some $j\in \{1, \ldots, N-1 \} $. In view of Lemmas \ref{primitive-of-p-over-q} and \ref{theta-equation}, the characteristic curves of $L$, in $S_j,$ can be written in the form:
\begin{equation}
  x(t,y) = y + c_j\ln |t-t_j|  + c_{j+1} \ln |t-t_{j+1}| + R_j(t).
  \label{chars2j}
\end{equation}
with $R_j\in C^\infty ([t_j,t_{j+1}])$ and $c_jc_{j+1} >0$.
We observe that $c_jc_{j+1}>0$ leads to
\begin{equation}
  \lim_{t\to t_j^+} x(t,y) =\lim_{t\to t_{j+1}^-}x(t,y) = \sign (c_j) \infty, \qquad y\in \R.
  \label{chars3j}
\end{equation}
Clearly \eqref{chars3j} implies that every smooth curve with endpoints on $t=t_j$ and $t=t_{j+1}$ is hit at least twice by the charateristic curve \eqref{chars2j} provided $y\gg 1$ (respectively, $-y\gg 1$) if $c_j>0$ (respectively, $c_j <0$), and therefore, the condition of Duistermaat-H\"ormander for the surjectivity fails.

Suppose now that there are no separatrix strips. Hence, $p(t)$ and $q(t)$ do not change sign in $[t_j,t_{j+1}]$, $j=0,1,\ldots, N$, $t_0=-\infty$, $t_{N+1} \doteq +\infty$ and  fixing $j$, we note that the line segment $x+\nu t = C$, $t\in [t_j, t_{j+1}]$ is transversal to $L$ provided $\nu\neq 0$ has the same sign as $p(t) q(t)$ for some $t\in ]t_j, t_{j+1}[$. So we have global picewise smooth global transversal. Smoothing by mollifiers $\veps^{-1} \varphi (\veps^{-1} t)$ near $t=t_j$   makes the curve smooth and still globally transversal provided $0< \veps \ll 1$. The proof of Theorem 1.1 is complete.

\vspace{10pt}

\begin{example}
\emph{We focus on the vector fields $L_{\lambda, k}$  defined in \eqref{vf5a} and exhibit some geometric features. The integral trajectories of $L_{\lambda, k}$ are given by the curves
\begin{equation}
  x(t) = \lambda\left[(-1)^k\frac{1}{2}\log\left|1+t\right| - \frac{1}{2}\log\left|1-t\right| - \sum_{i<k}^{\rule{4mm}{0.1mm}} \frac{t^i}{i}\right]
\end{equation}
where $\overline{\sum}$ extends only to odd numbers when $k$ is even and only to even numbers
when $k$ is odd.}

{\em The vector fields $L_{\lambda, k}$ are intrinsically Hamiltonian vector fields, i.e.
they are tangent to the level sets of a regular smooth function on the plane -- equivalently, the kernel of each operator $L_{\lambda, k}$ contains regular smooth functions.

For example, the following smooth function $f_{\lambda, k} \in\ker (L_{\lambda, k})$:
\begin{eqnarray}
  f_{\lambda, 2k+1}(x,t) &=& (1-t^2)\exp\left[2\left(\frac{x}{\lambda}+\sum_{i<2k+1}^{\rule{4mm}{0.1mm}} \frac{t^i}{i}\right)\right], \mbox{ and }\\
  f_{\lambda, 2k}(x,t) &=&\tan^{-1}\left\{\displaystyle\frac{1-t}{1+t} \exp\left[2\left( \frac{x}{\lambda}+ \sum_{i<2k}^{\rule{4mm}{0.1mm}} \frac{t^i}{i}\right)\right]\right\}.
\end{eqnarray}
}
%In Fig.~\ref{fig:k1k2} we show the level sets for $L_{1,1}$ and $L_{1,2}$. As the pictures suggest,
%the $L_{\lambda, k}$ admit no global transversal for $k$ odd while they do if $k$ is even.
%Indeed, let $G_{\lambda, k}:\R^2\to\R P^1$ be the Gauss map associating to each point the
%direction of $L_{\lambda, k}$ in that point.

\end{example}
\begin{figure}
  \begin{center}
  \includegraphics[width=5cm]{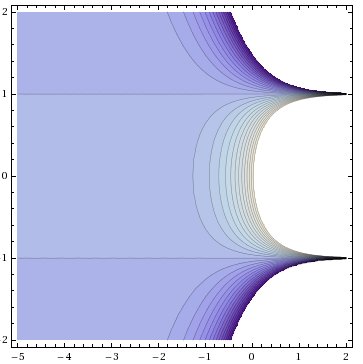} \hskip.75cm \includegraphics[width=5cm]{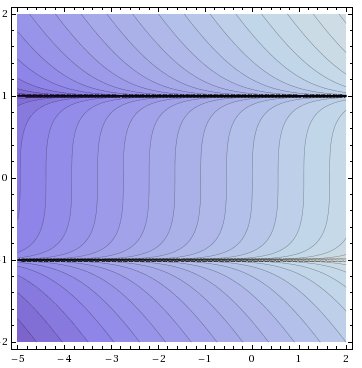}
  \end{center}
   { \caption{Integral curves of $L_{1,1}=(1-t^2)\partial_t+t\partial_x$ and $L_{1,2}=(1-t^2)\partial_t+t^2\partial_x$, respectively. Clearly no global transversal exists for $L_{1,1}$ while $L_{1,2}$ is topologically equivalent to a constant vector field.}}
  \label{fig:k1k2}
\end{figure}

\begin{remark} \emph{We can generalize Theorem 1.2 for smooth non-singular vector fields assuming that $p$ and $q$ are in general position with respect to each other, i.e., each zero of $p$ and $q$ has finite multiplicity.
Choosing $t_1$ and $t_2$ to be two successive zeros of $p(t),$ then $t_1$ and $t_2$ form a separatrix if and only if the sum of degrees of all the roots of $q$ between $t_1$ and $t_{2}$ is odd.}
\end{remark}

\section{Estimates on the right inverse}\

The aim of this section is to prove the Theorem \ref{Theorem2}. First we will construct a right inverse as follows:

Let $j\in \{1, \ldots, N-1 \}$. If  the strip  $S_j$ is a separatrix, we use Lemma \ref{primitive-of-p-over-q} to obtain
\begin{eqnarray}
L^{-1}_j f & = & \int_{\theta_j}^t \frac{f(\tau, x + \rho (\tau) -\rho (t))}{p(\tau)}  \, d\tau\nonumber\\
& = & \int_{\theta_j}^t \frac{f(\tau, c_j \ln \frac{|\tau-t_j|}{|t-t_j|} + c_{j+1} \ln \frac{|\tau-t_{j+1}|}{|t-t_{j+1}|}+ R_j(\tau)-R_j(t))}{p(\tau)}  \, d\tau
\label{invr1}
\end{eqnarray}

If $S_j$ is not separatrix, we construct $L^{-1}_j $ as  the Green function for the Cauchy problem in $S_j$
\begin{eqnarray}
L^{-1}_j f (t,x) & = & G_j^\nu f(t,x),
\label{invr1G}
\end{eqnarray}
where $\nu \neq 0$ is fixed by the requirement $C_j: x+\nu t =0$, $t\in [t_j, t_{j+1}]$ is noncharacteristic for $L$ in $S_j$ and  $u_j (t,x) = G_j^\nu f(t,x) $ si defined by
\begin{equation}
Lu_j = f, \ \ (t,x) \in S_j, \quad u|_{C_j} =0.
\label{invr2G}
\end{equation}

The global transversality of $C_j$ in $S_j$ implies that $u_j\in C^\infty (\, \overline{S}_j )$ (we are in a particular case of \cite{DH1}).

The next assertion plays a crucial role in the proof of the global solvability for $L$ in the presence of the separatrix strip.\\

\begin{proposition}
Suppose that  $S_j$ is a  separatrix and set $I_j \doteq ]t_j,t_{j+1}[$, then $L^{-1}_j $ has the following properties:

\begin{enumerate}
  \item[{\it i)}] If $C^\infty (\, I_j: E^\veps_{gr} (\R))$ (respectively, $C^\infty (\, I_j: E^\veps_{dec} (\R))$) is the subspace of $C^\infty(I_j~\times~\R)$ consisting of all infinitely differentiable functions that satisfy the following growth (respectively, decay) condition
      \begin{equation}\label{exponential-growth-veps}
        \forall \alpha \in \Z_+^2, \exists C>0 \textit{ such that } |\partial_{t,x}^\alpha f(t,x) |\leq Ce^{\veps |x|}, \, t\in I_j, x\in \R
      \end{equation}
      (respectively,
       \begin{equation}\label{exponential-decay-veps}
        \forall \alpha \in \Z_+^2, \exists C>0 \textit{ such that } |\partial_{t,x}^\alpha f(t,x) |\leq Ce^{-\veps |x|}, \, t\in I_j, x\in \R )
      \end{equation}
  then
      \begin{equation}
        L^{-1}_j: C^\infty (\, \overline{I}_j: E^\veps_{gr} (\R)) \longrightarrow   L^1 (I_j: E^\veps_{gr} (\R))\, \mbox{$\bigcap$}\,  C^\infty ( I_j :E^\veps_{gr} (\R)) \label{invr4j}
      \end{equation}
      (respectively,
      \begin{equation}
        L^{-1}_j: C^\infty (\, \overline{I}_j: E^\veps_{dec} (\R)) \longrightarrow   L^1 (I_j: E^\veps_{dec} (\R))\, \mbox{$\bigcap$}\,  C^\infty ( I_j :E^\veps_{dec} (\R)) \label{invr5j})
      \end{equation}
      if  
      \begin{equation}\label{eps-1}
0<\veps < \min \{|c_j|^{-1}, |c_{j+1}|^{-1} \}
      \end{equation}

  \item[{\it ii)}] For $s_1,s_2\in \R,$
      \begin{equation}
        L^{-1}_j: C(\, \overline{I}_j: H^{s_1,s_2} (\R)) \longrightarrow   L^1 (I_j:   H^{s_1,s_2} (\R))\, \mbox{$\bigcap$}\, C(I_j: H^{s_1,s_2}(\R)) \label{invr6j}
      \end{equation}
  Moreover, for any $\veps >0$ we have
      \begin{equation}\label{sup-tj}
        \sup_{t \in [t_j,t_{j+1}]}\left(|t-t_j|^\veps |t-t_{j+1}|^\veps \| L_j^{-1} f(t,\cdot) \|_{H^{s_1,s_2} (\R)}  \right) \leq C_{\veps, s_1, s_2} \| f \|_{ C(\, \overline{I}_j: H^{s_1,s_2} (\R))}
      \end{equation}

  \item[{\it iii)}] If $ f(t,x) = \sum_{\ell=0}^k f_\ell(t) x^\ell$, then
      \begin{eqnarray}
        L^{-1}_j f(t,x) = \sum_{\ell=0}^k g_\ell(t) x^\ell\label{inv7j}
      \end{eqnarray}
  with
      \begin{align}
        g_k(t) & =  g_k(t_\mu) \gamma_\mu \ln^{k+1}\mbox{$\frac{1}{|t-t_\mu|}$} (1 + o(1)) \textrm{ near $t= t_\mu$, $\gamma_\mu \neq 0$, $\mu =j, j+1$, } \label{invr8j}\\
       g_\ell(t)& = o(\ln^{k+1} \mbox{$\frac{1}{|t-t_\mu|}$}) \textrm{ near $t= t_\mu$, $\mu = j, j+1$}, \ \ \ell =0,1,\ldots, k-1. \label{invr8j} \\
      \end{align}

  \item[{\it iv)}] Given a zero order p.d.o. $b(t,x,D)$ in $x$, smoothly depending on $t$, and $s_1,s_2\in \R,$ we can find $\veps_0 = \veps_0(L,s_1,s_2) >0$ such that if
      \begin{equation}
        \max_{{|\alpha|\leq [s_1]+2}\atop{|\beta| \leq [s_2]+2}}\sup_{{t\in [t_j,t_{j+1}]}\atop{(x,\xi)\in \R^2}} \langle x \rangle^{-\alpha} \langle \xi \rangle^{-\beta} |\partial_x^\alpha\partial_\xi^\beta b(t,x,\xi)|   < \veps_0 \label{pdo1}
      \end{equation}
  then $L+b(t,x,D)$ admits a right inverse which satisfies \eqref{invr6j}. \\
\end{enumerate}
\end{proposition}

\begin{proof}
We observe that for $t$ close to $t_j$ we can write
\begin{eqnarray}
L^{-1}_j f (t,x) & = & \int_{\theta_j}^t \, f(\tau, c_j \ln \frac{|\tau-t_j|}{|t-t_j|} + M_j (t,\tau) )
\frac{\sigma_j(\tau) }{\tau -t_j}  \, d\tau
\label{invrest1}
\end{eqnarray}
with $\sigma_j \in C^\infty ([t_j, \theta_j])$, $M_j\in C^\infty (\Delta_j)$, $\Delta_j = \{ t_j \leq \tau \leq t\leq \theta_j \}$.
Therefore,
\begin{eqnarray}
|\partial_x^\alpha L^{-1}_j f (t,x)| & \leq  & \int^{\theta_j}_t \, |\partial^\alpha_x f(\tau, x+ c_j \ln \frac{|\tau-t_j|}{|t-t_j|} + M_j (t,\tau) )
\frac{\sigma_j(\tau) }{\tau -t_j} | \, d\tau\nonumber\\
& \leq  & C e^{\veps |x|} \int^{\theta_j}_t \, e^{\veps |c_j| \ln \frac{\tau-t_j}{t-t_j} }
\frac{1}{\tau -t_j}  \, d\tau\nonumber\\
& =  & C e^{\veps |x|} \frac{1}{(t-t_j)^{\veps |c_j|} }\int^{\theta_j}_t \,
\frac{1}{(\tau -t_j)^{1-\veps|c_j|}}  \, d\tau\nonumber\\
& =  & C e^{\veps |x|} \frac{1}{\veps |c_j| (t-t_j)^{\veps |c_j|} } ((\theta_j  -t_j)^{\veps|c_j|}  - (t  -t_j)^{\veps|c_j|})
\nonumber\\
 & =  & C e^{\veps |x|} \frac{(\theta_j  -t_j)^{\veps|c_j|} }{\veps |c_j| (t-t_j)^{\veps |c_j|} } (1  + O( (t  -t_j)^{\veps|c_j|}))
\label{invrest2}
\end{eqnarray}
Similarly, we derive that near $t_{j+1}$ we have
\begin{equation}
|\partial_x^\alpha L^{-1}_j f (t,x)| \leq
  C e^{\veps |x|} \frac{(t_{j+1}-\theta_j )^{\veps|c_{j+1}|} }{\veps |c_{j+1}| (t_{j+1}-t)^{\veps |c_{j+1}|} } (1  + O( (t_{j+1}-t)^{\veps|c_{j+1}|}))
\label{invrest3}
\end{equation}
Clearly, \eqref{invrest1}, \eqref{invrest2}, \eqref{invrest3} imply  \eqref{invr4j} provided $0<\veps < \min \{\frac{1}{|c_j|}, \frac{1}{|c_{j+1}|} \}$.\\

\noindent As it concerns to item {\it ii)}, taking into account the inequality
\begin{eqnarray}
\sup_{x\in \R, |\lambda| \geq 1}  |\lambda|^{-|s_2|} \langle x\rangle^{s_2}  \langle x+ \lambda \rangle^{-s_2} & < & +\infty
\label{ineq1}
\end{eqnarray}
we observe that for $\alpha\in \Z_+$ and $s_2\in \R$ we have for $t$ near $t_j$
\begin{eqnarray}
\|\langle \cdot \rangle^{s_2} \partial_x^\alpha L^{-1}_j f (t,\cdot )\|_{L^2} & \leq  & C \int^{\theta_j}_t \, \sup_{x\in R} \left (\langle x\rangle ^{s_2}\langle x+ c_j \mbox{$\ln \frac{|\tau-t_j|}{|t-t_j|}$} \rangle^{-s_2} \right) \mbox{$\frac{1}{|\tau -t_j|}$}  \, d\tau \nonumber\\
 &  & \times \, \sup_{t \in [t_j,t_{j+1}]}\! \| \langle \cdot \rangle^{s_2}\partial^\alpha f(t, \cdot) \|_{L^2} \nonumber \\
 & \leq  & \tilde{C}\int^t_{\theta_j} \ln^{|s_2|}\mbox{$\left( \frac{\tau - t_j}{t-t_j} \right)\frac{1}{\tau - t_j}$}d\tau \!\! \sup_{t \in [t_j,t_{j+1}]}\! \| \langle \cdot \rangle^{s_2}\partial^\alpha f(t, \cdot) \|_{L^2} \nonumber\\
  & =  & \frac{\tilde{C}}{|s_2|}\ln^{|s_2|+1}\frac{1}{t-t_j} \sup_{t \in [t_j,t_{j+1}]} \| \langle \cdot \rangle^{s_2}\partial^\alpha f(t, \cdot) \|_{L^2}
\label{invrest4}
\end{eqnarray}
Therefore we obtained \eqref{invr6j} for $s_1 \in \Z_+$ (summation in \eqref{invrest4} over $|\alpha|$). We conclude the general case for $s_1$ by interpolation and duality arguments.

Since the logarithmic singularity is weaker then any polynomial one, \eqref{invrest4} yields \eqref{sup-tj}
\end{proof} \

Next, we show a gluing lemma, which will imply  that
\begin{eqnarray}
L^{-1} f (t,x) & = & L_j^{-1} f(t,x), \qquad (t,x) \in S_j, \, j=0,1, \ldots, N
\label{invr5}
\end{eqnarray}
is a right inverse satisfying the properties stated in Theorem 1.2. This gluing auxiliary assertion seems to be also  a novelty ``per se" and might be of an independent interest.

Let $\Omega $ be an open domain in $\R^n$ and let $\delta>0.$
Set $I_\delta = ]-\delta, \delta[$, $I^+_\delta = ]0,\delta[$, $I^-_\delta = ]-\delta, 0[$, and
\begin{eqnarray}
\Omega^\pm_\delta = I^\pm_\delta \times \Omega =\{ (t,x):\,  0<\pm t < \delta,  x\in \Omega \},
\label{dom1}\\
\Omega_\delta = I_\delta \times \Omega = \{ (t,x):\,  |t| < \delta,  x\in \Omega \}.
\label{dom2}
\end{eqnarray}

Consider the smooth vector field
\begin{eqnarray}
X = a(t,x) \partial_t + \sum_{j=1}^n a_j(t,x) \partial_{x_j},
\label{op1}
\end{eqnarray}
having $t=0$ as a characteristic, i.e.,
\begin{equation}
a_0(0,x) = 0, \quad x\in \Omega
\label{op2}
\end{equation}

Let
\begin{equation}
b= b(t,x) \in C^\infty (\Omega_\delta)
\label{op2}
\end{equation}
or, in the case $\Omega = \R^n$ we allow
\begin{equation}
\textrm{$b$ to be a zero order p.d.o. in $x$ (cf. \cite{Co1}) depeding smoothly on $t\in ]-\delta, \delta[$}
\label{op3}
\end{equation}

We have

\begin{lemma}\label{thm:AbstractLemma}
Let $f \in C^\infty (\Omega_\delta)$ (respectively, $f\in C(\, ]-\delta, \delta[\, : H^{s_1,s_2}(\R^n))$ for some $s_1,s_2\in \R$ if $\Omega = \R^n$).
Suppose that
\begin{eqnarray}
u^\pm &\in & C^\infty  (\Omega^\pm_\delta)
\label{wsola}
\end{eqnarray}
(respectively,
\begin{eqnarray}
u^\pm &\in & C (I^\pm_\delta: H^{s_1,s_2} (\R^n))
\label{wsolb}
\end{eqnarray}
for some $s_1,s_2\in \R$)
 satisfies
\begin{eqnarray}
X u^\pm + bu^\pm  = f \ \ \ \textrm{in $ \Omega^\pm_\delta$}
\label{wsol1}
\end{eqnarray}

Then
\begin{eqnarray}
u(t,x) = \left\{ \begin{array}{ccc} u^+(t,x)  & \textrm{if} & (t,x)\in  \Omega^+_\delta \\
u^-(t,x)  & \textrm{if} & (t,x)\in  \Omega^-_\delta
\end{array} \right.
\label{wsol2}
\end{eqnarray}
is  a well defined $L^1_{loc} (\Omega)$ (respectively, $L^1 (I_\delta: H^{s_1,s_2}(\R^n) )$
distributional solution of $Xu=f$ in $\Omega_\delta$ provided
\begin{eqnarray}
u^\pm &\in & L^1(I^\pm_\delta\times K), \qquad K\subset\subset \Omega
\label{wsol3a}
\end{eqnarray}
(respectively,
\begin{eqnarray}
u^\pm &\in & L^1(I^\pm_\delta:H^{s_1,s_2}(\R^n)),
\label{wsol3b}
\end{eqnarray}
if $\Omega = \R^n)$)
and
\begin{eqnarray}
  \lim_{t\to 0^\pm }\int_{\R^n} a(t, x) u^\pm (t,x)\varphi(t,x) dx   = 0, \qquad \varphi \in C^\infty_0(\Omega_\delta)
  \label{wsol4}
\end{eqnarray}

\end{lemma}
\begin{proof}
Let $\varphi(t,x) \in C_0^\infty (\Omega_\delta)$. We have to prove that
\begin{eqnarray}
\langle u, X^* \varphi + b^* \varphi \rangle    = \langle f, \varphi \rangle
  \label{wsol5}
\end{eqnarray}
where $X^*$  (respectively, $b^*$) stands for the adjoint of $X$ (respectively, $b$).

Taking into account \eqref{wsol3a}, \eqref{wsol3b} and Lebesgue's dominated convergence theorem we have
\begin{equation}\label{wsol5a}
\langle u, X^* \varphi + b^* \varphi \rangle    = \lim_{\veps \to 0} (J^+_\veps(u^+,\varphi)+J^-_\veps(u^-,\varphi)),
\end{equation}
where
\begin{equation}
J^\pm_\veps (u^\pm ,\varphi) = \pm \int_{\pm \veps}^{\pm \delta} \left(
\int_\Omega u^\pm (t,x) (X^*\varphi (t,x)  + b^*(t,x,D) \varphi (t,x)) dx \right) dt.
    \label{wsol7}
\end{equation}

Integration by parts, duality arguments, the Fubini theorem  and  \eqref{wsol1} imply that
\begin{eqnarray}
J^\pm_\veps (u^\pm ,\varphi)    &=& \pm \int_{\pm \veps}^{\pm \delta}
\int_\Omega (X u^\pm (t,x) +  b(t,x,D)u^\pm (t,x)) \varphi (t,x) dx  dt \nonumber\\
&  & + \int_\Omega a(\pm \veps, x) u^\pm (\pm \veps, x) \varphi (\pm \veps, x) dx \nonumber\\
& = & \int_{\Omega^\pm_\delta \setminus \overline{\Omega^\pm_\veps}} f(t,x) \varphi (t,x) \nonumber\\
&  & + \int_\Omega a(\pm \veps, x) u^\pm (\pm \veps, x) \varphi (\pm \veps, x) dx
    \label{wsol7}
\end{eqnarray}
Next, using the hypothesis \eqref{wsol4}, we deduce that
\begin{eqnarray}
\lim_{\veps \to 0}J^\pm_\veps (u^\pm ,\varphi)   & = & \int_{\Omega^\pm_\delta } f(t,x) \varphi (t,x) dt dx
    \label{wsol8}
\end{eqnarray}
and, plugging into the RHS of \eqref{wsol5},  we obtain,
\begin{eqnarray}
\langle u, L^* \varphi + b^* \varphi \rangle    &=& \int_{\Omega^+_\delta } f(t,x) \varphi (t,x) dtdx + \int_{\Omega^-_\delta } f(t,x) \varphi (t,x) dtdx \nonumber\\
& = & \int_{\Omega_\delta } f(t,x) \varphi (t,x) dtdx
     \label{wsol9}
\end{eqnarray}
This completes the proof of the lemma.
\end{proof}

Combining Proposition 3.1 and Lemma 3.2 we derive the assertions for $L^{-1}$.

As it concerns the perturbation with $b(t,x,D)$, we reduce the equation in $\R^2$ to $Lu + b(t,x,D_x)u= f$ on $S_j$, $j=0,1,\ldots,N $. We are reduced to the study of the global solvability of
 \begin{eqnarray}
u +   L^{-1}b(t,x,D)u = L^{-1}f, \qquad (t,x) \in S_j, j=0,1,\ldots, N.
     \label{decomp1}
\end{eqnarray}
We apply the  Picard type scheme
\begin{eqnarray}
u_k =-   L^{-1}b(t,x,D)u_{k-1} + L^{-1}f,\qquad  k\in \N, u_0 = 0
     \label{decomp1a}
\end{eqnarray}

If $j=1, \ldots, N$, we use the results for $H^{s_1,s_2}$ estimates of p.d.o. in $\R^n$ (e.g., cf. \cite{Co1}) and choose $\veps_0$ so small that
\begin{eqnarray}
\| b(t,x,D)  L^{-1}\|_{L^1([t_1,t_N]: H^{s_1,s_2})\to L^1([t_1,t_N]:H^{s_1,s_2})} < 1
     \label{decomp2}
\end{eqnarray}
Using continuity arguments we can find $\delta >0$ (small enough) such that
\begin{eqnarray}
\|L^{-1} b(t,x,D)  L^{-1}\|_{L^1([t_-\delta,t_N+\delta]: H^{s_1,s_2})\to L^1([t_1-\delta,t_N+\delta]:H^{s_1,s_2})} < 1
     \label{decomp2a}
\end{eqnarray}
Since $p(t)$ has no zeroes for $t>t_N+\delta$ and $t\leq t_1-\delta$ we have the following estimates: there exista $C=C_\delta >0$ such that
\begin{eqnarray}
\|L^{-1} bu (t,\cdot)\|_{H^{s_1,s_2}} \leq C_\delta \int_{\theta_j}^t \| u(\tau, \cdot) \|_{H^{s_1,s_2}} d\tau,
     \label{decomp3}
\end{eqnarray}
for $j=0$, $t\leq t_1-\delta$, $j=N$, $t\geq t_N+\delta$. Combination of contraction and Gronwall inequlaities (see \cite{GG1}) imply the convergence of \eqref{decomp1a} and the existence of $(L+b)^{-1}$ satisfying the last part of Theorem 1.2

\begin{remark}
\emph{We point out that the estimates for $f \in C^\infty (\, \overline{I}_j: E^\veps_{dec} (\R))$ allows to extend solvability for $L$ and $L+b$ in Gelfand-Shilov spaces $S^\mu_\mu(\R)$ in $x$, provided $\mu > 1.$ See \cite{CGR1} for global solvability and regularity results for some degenerate p.d.o. under similar subexponential decay conditions. We can show that, if the decay is superexponential the solution $u$ loses this decay, unlike the solvability in Gelfand-Shilov spaces $S^\mu_\mu, 1/2 \leq \mu \leq 1,$ cf. see \cite{DW1,GPR1} and the references therein.}
\end{remark}

\section{The sharpness of the estimates for $L_0$}\

We consider the model equation $L_0u = f$. Using the method of the characteristics, for $t \neq \pm 1,$ one can write formally a right inverse of $L_0$ in the following way
\begin{equation}
L^{-1}_0 f \doteq  \int_{0}^t \, f(\tau ,x  +  \ln |\frac{1-\tau^2}{1-t^2}|)  \frac{1}{1-\tau^2} \, d\tau
 =  G_+ f + G_-f,
\label{modvf3}
\end{equation}
 where
\begin{equation}
G_\pm f(t,x) \doteq \frac{1}{2}\int_{0}^t \, f(\tau, x  +  \ln |\frac{1-\tau}{1-t} | + \ln |\frac{1+\tau}{1+t}| )  \frac{1}{1 \pm \tau} \, d\tau
\label{modvf4}
\end{equation}\

We define in a natural way $C^\infty (\R: E^\veps_{gr} (\R))$ as the inductive limit 
\begin{equation}
C^\infty (\R: E^\veps_{gr} (\R))  = \lim_{T\nearrow +\infty} C^\infty (\, [T,T]: E^\veps_{gr} (\R))
\label{modvf5}
\end{equation}

Observe that $C^\infty (\R:E^\veps_{gr}(\R))$ is a vector subspace of $C^\infty (\R^2)$ and, given $f_1, f_2 \in C^\infty (\R: E^\veps_{gr} (\R)),$ we have $f_1\!\cdot\!\! f_2 \in  C^\infty (\R:E^\veps_{gr}(\R)).$ In particular, the projections $\pi_1(t,x) = t$ and $\pi_2(t,x) = x$ \ belong to this space and consequently, any polynomial function $p$ belongs to $C^\infty (\R: E^\veps_{gr} (\R)).$

We introduce a topology on $C^\infty (\R: E^\veps_{gr} (\R))$ by the following family of seminorms
\begin{equation}\label{seminorms}
    \rho_{j,k,T}^{\, \veps}(f) \doteq \sup \{\ |e^{-\veps |x|}\partial^{\alpha_1}_{t}\partial^{\alpha_2}_{x} f(t,x)|;  |\alpha_1| \leq j, |\alpha_2| \leq k, |t| \leq T, x \in \R, \} \\
\end{equation}
where $T>0$ and $j,k \in \Z_+.$

\begin{lemma}\label{estimate_int}
If $a \in C^1(\R)$ and $p \in \N$ then, when $t \rightarrow 1,$ we have
$$
\int_{0}^t a(s)\ln^p \left| \frac{1-s}{1-t}\right| \frac{1}{1-s}\, ds = \frac{a(1)}{p+1} \ln^{p+1} \left| \frac{1-s}{1-t}\right|(1+o(1))
$$
\end{lemma}
\begin{proof}
\begin{eqnarray*}
  \int_{0}^t a(s)\ln^p \mbox{$\left| \frac{1-s}{1-t}\right|\frac{1}{1-s}$}\, ds &=& a(1)\int_{0}^t \ln^p \mbox{$\left| \frac{1-s}{1-t}\right|\frac{1}{1-s}$}\ ds + \int_{0}^t a_1(s)\ln^p \mbox{$\left| \frac{1-s}{1-t}\right|$} ds \\
  &=& \frac{a(1)}{p+1}\ln^{p+1} \mbox{$\left| \frac{1-s}{1-t}\right|$} + o\left(\ln^p \mbox{$\frac{1}{|1-t|}$} \ \right) \\
  &=& \frac{a(1)}{p+1}\ln^{p+1} \mbox{$\left|\frac{1}{1-t}\right|$}(1+o(1))
\end{eqnarray*}\vskip-2em
\end{proof}

\begin{lemma}\label{monomial}
If f is a monomial function with respect to $x$, i.e., $f(t,x) = f_j(t) x^j,$ with $f_j \in C^1(\R)$ and $j \in \Z_+,$ then \begin{equation}\label{inverse_monomial}
L^{-1}_0 f(t,x) = \sum_{\ell=0}^j g_{j\ell}(t) x^\ell
\end{equation}
with
\begin{align}
g_{j0}(t) & = \frac{f_0(\pm 1)}{2} \ln \left|\frac{1}{1\mp t} \right| (1+o(1)), && t\to \pm 1\\
g_{j\ell} (t)& = O( \ \ln^{j+1-\ell} \frac{1}{|1 \mp t|} \ ),  && t\to \pm 1.
\end{align}
\end{lemma}
\begin{proof}
From (\ref{modvf3}) and (\ref{modvf4}) we obtain
\begin{eqnarray*}
  G_\pm f(t,x) & = & \frac{1}{2}\int_{0}^t \, f_j(\tau) \left( x  +  \ln \left|\frac{1-\tau}{1-t}\right| + \ln \left|\frac{1+\tau}{1+t}\right| \right)^j \frac{1}{1 \pm \tau} \, d\tau \\
  & = & \sum_{\ell=0}^{j}\left[\frac{1}{2}\binom{j}{\ell} \int_{0}^t \, f_j(\tau) \left(\ln \left|\frac{1-\tau}{1-t}\right| + \ln \left|\frac{1+\tau}{1+t}\right| \right)^{j-\ell}  \frac{1}{1 \pm \tau} \, d\tau \right]x^{\ell}  \\
   & = & \sum_{\ell=0}^{j}g_{{j\ell}\pm}(t)\, x^{\ell}  \\
\end{eqnarray*}
where
\begin{eqnarray*}
   g_{{j\ell}\pm}(t) & = & \frac{1}{2}\binom{j}{\ell} \int_{0}^t \, f_j(\tau) \left(\ln \left|\frac{1-\tau}{1-t}\right| + \ln \left|\frac{1+\tau}{1+t}\right| \right)^{j-\ell}  \frac{1}{1 \pm \tau} \, d\tau \nonumber \\
   & = & \frac{1}{2}\binom{j}{\ell}\sum_{m=0}^{j-\ell}\binom{j-\ell}{m}  \int_{0}^t \, f_j(\tau) \ln^m \left|\frac{1-\tau}{1-t} \right| \ln^{j-\ell-m} \left|\frac{1+\tau}{1+t}\right| \frac{1}{1 \pm \tau} \, d\tau
\end{eqnarray*}

Now, it follows from Lemma \ref{estimate_int} that, near $t=1$, we have

\begin{eqnarray*}
   g_{{j\ell}-}(t) & = &\frac{1}{2}\binom{j}{\ell}\sum_{m=0}^{\ell}\binom{j-\ell}{m}  \int_{0}^t \, f_j(\tau) \ln^m \left|\frac{1-\tau}{1-t} \right| \ln^{j-\ell-m} \left|\frac{1+\tau}{1+t}\right| \frac{1}{1 - \tau} \, d\tau \\
   & = & \frac{1}{2}\binom{j}{\ell}\sum_{m=0}^{j-\ell}\binom{j-\ell}{m} \frac{f_j(1)\ln^{j-\ell-m} \left|\frac{2}{1+t}\right|}{m+1} \ln^{m+1} \left|\frac{1}{1-t} \right| (1+o(1))\\
   & = & O\left( \ \ln^{j+1-\ell} \left|\frac{1}{1-t} \right| \ \right)
\end{eqnarray*}
Analogously, near $t=-1$, we have
\begin{eqnarray*}
   g_{{j\ell}+}(t) & = &\frac{1}{2}\binom{j}{\ell}\sum_{m=0}^{j-\ell}\binom{j-\ell}{m}  \int_{0}^t \, f_j(\tau) \ln^m \left|\frac{1-\tau}{1-t} \right| \ln^{j-\ell-m} \left|\frac{1+\tau}{1+t}\right| \frac{1}{1 + \tau} \, d\tau \\
   & = & \frac{1}{2}\binom{j}{\ell}\sum_{m=0}^{j-\ell}\binom{l}{m} \frac{f_j(1)\ln^{m}|\frac{2}{1-t}|}{j-\ell-m+1} \ln^{j-\ell-m+1} \left|\frac{1}{1+t}\right|(1+o(1))\\
   & = & O\left( \ \ln^{j+1-\ell} \left|\frac{1}{1+t} \right| \ \right)
\end{eqnarray*}
In particular, for $\ell=0,$ we have
\begin{equation*}\label{g0-pm}
   g_{j0\pm}(t) = \frac{1}{2} \int_{0}^t \, f_j(\tau) \frac{1}{1 \pm \tau} \, d\tau= \frac{f_j(\pm 1)}{2} \ln \left|\frac{1}{1\mp t} \right| (1+o(1)).
\end{equation*}\vskip-1.5em
\end{proof}

The next assertion shows that we have sharp estimates on the singularities.

\begin{proposition}
  The following properties hold: there exists $\veps_0>0$ such that for all $\veps \in ]0,\veps_0[$
\begin{enumerate}
 \item[i)] Given $T>0$ and $k \in \Z_+$, we have
  \begin{equation}
 \sup_{{x \in \R}\atop{|\alpha|\leq k}} |(1-t^2)^\veps L^{-1}_0f(t,\cdot)| \leq \, C^{\veps_0}_{T} \, \rho^{\veps_0}_{0,k,T}(f)
  \end{equation}

 \item[ii)] If $ f(t,x) = \sum_{j=0}^k f_j(t) x^j
   \label{modvf8}$, then
  \begin{equation}
   L^{-1}_0 f(t,x) = \sum_{j=0}^k g_j(t) x^j\label{modvf9}
  \end{equation}
 with
  \begin{align*}
   g_0(t) &  =  \frac{f_k(\pm 1)}{2(k+1)} \ln^{k+1} |1\pm t| (1 + o(1)), & t\to \pm 1, \label{modvf10}\\
   g_j (t)& =  O ( \ln^{k+1-j} |1\pm t|), & t\to \pm 1.
  \end{align*}
 \item[iii)] $u=L^{-1}_0 f$ is a weak solution of $Lu = f$ for all $f\in C^\infty (\R: E^\veps (\R))$ such that $\forall \alpha \in \Z_+$, $K\subset\subset \R$,   there exist $M>0$ such that
  \begin{equation}
    |\partial_x^\alpha u(t,x)| \leq M |1\pm t|^{-\veps}, \quad 0 < |1\pm t| \ll 1,\, x\in K
    \label{mod1}
  \end{equation}
\end{enumerate}
\end{proposition}

%==================================================== beginning of the proof ==================
\begin{proof}

To prove \emph{i)} we start by defining, for each $T>0, k \in \Z_+$ and $u \in C^\infty (\R:E^\veps(\R))$ the following function:
$$
P^{\veps_0}_{k,T}(u) \doteq \int_{-T}^{T} \ \sup_{{x \in \R}\atop{|\alpha|\leq k}}\left| e^{-\veps_0 |x|}\partial^{\alpha}_{x} u(t,x) \right| dt
$$
Thus, for any $f \in C^\infty (\R:E^\veps(\R)),$ with $0 < \veps < \veps_0$ and $t>0$ we have
\begin{eqnarray*}
 & & P^{\veps_0}_{k,T}(G_- f) =  \int_{-T}^{T} \ \sup_{{x \in \R,}\atop{|\alpha|\leq k}}\left| e^{-\veps_0 |x|}\partial^{\alpha}_{x} G_- f(t,x) \right| dt \\
  & = & \int_{-T}^{T} \ \sup_{{x \in \R,}\atop{|\alpha|\leq k}}\left| e^{-\veps_0 |x|}\partial^{\alpha}_{x}\left[ \frac{1}{2}\int_{0}^t \, f(\tau, \mbox{$x+\ln |\frac{1-\tau}{1-t}|+\ln |\frac{1+\tau}{1+t}|$} )  \frac{1}{1 - \tau} \, d\tau \right] \right| dt \\
    & \leq & \frac{1}{2} \int_{-T}^{T} \, \int_{0}^t \sup_{{x \in \R,}\atop{|\alpha|\leq k}}\left| e^{-\veps_0 |x|}\partial^{\alpha}_{x} f(\tau, \mbox{$x+\ln |\frac{1-\tau}{1-t}|+\ln |\frac{1+\tau}{1+t}|$} )  \frac{1}{1 - \tau} \, d\tau \right| dt \\
  & \leq & \frac{1}{2} \rho^{\veps}_{0,k,T}(f) \int_{-T}^{T}\int_{0}^t\sup_{{x\in\R,}\atop{|\alpha|\leq k}}\left| e^{-\veps_0|x|}\exp(\veps(\mbox{$|x|+\ln |\frac{1-\tau}{1-t}|-\ln |\frac{1+\tau}{1+t}|$ })) \frac{1}{1 - \tau} \right| d\tau dt \\
    & \leq & \frac{1}{2}\, e^{\veps-\veps_0} \rho^{\veps}_{0,k,T}(f) \int_{-T}^{T} \int_0^t \, \left|\frac{1+\tau}{1+t} \right|^{-\veps} \left|\frac{1-\tau}{1-t}\right|^\veps \frac{1}{|1-\tau|} d\tau dt \\
  & \leq & \frac{1}{2} \, C^{\veps_0}_{T} \, \rho^{\veps_0}_{0,k,T}(f)(1-t)^{-\veps}
\end{eqnarray*}

By using the same arguments, when $t<0,$ we obtain an analogous estimate to $G_+ f,$ and consequently
$$
P^{\veps_0}_{k,T}(L^{-1}_0f(t,x)) \leq C^{\veps_0}_{T} \rho^{\veps_0}_{0,k,T}(f)(1-t^2)^{-\veps}
$$

To prove \emph{ii)}, we use the results in the lemmas \ref{estimate_int} and \ref{monomial} below to obtain
\begin{eqnarray*}
  L^{-1}_0 f(t,x) &=&  \sum_{\ell =0}^{k} L^{-1}_0 (f_j(t)x^j)\ = \ \sum_{\ell =0}^{k} \left( \sum_{\ell =0}^j g_{j\ell}(t) x^\ell \right)\\
  &=& \sum_{\ell =0}^{k} \left( \sum_{\ell =j}^k g_{\ell j}(t) \right)  x^\ell \ = \ \sum_{\ell=0}^{k} g_{j}(t)  x^\ell
\end{eqnarray*}
where $$\displaystyle g_{j}(t) \doteq  \sum_{\ell =j}^k g_{\ell j}(t) = \sum_{\ell =j}^k O\left(\ln^{\ell +1-j}\frac{1}{|1\mp t|}\right) = O\left(\ln^{k+1-j}\frac{1}{|1\mp t|}\right), \mbox{ when } t \to \pm 1$$

To prove the statement \emph{iii)}, first, for $0<t<1,$ we have

\begin{eqnarray}
 |\partial^\alpha_x(G_- f(t,x))| & \leq &  \frac{1}{2}\int_{0}^t \, \left| \partial^{\alpha}_{x}f(\tau,
 \mbox{$x+\ln |\frac{1-\tau}{1-t}|+\ln |\frac{1+\tau}{1+t}|$})  \frac{1}{1 - \tau} \right| d\tau  \nonumber \\
   & \leq &  \frac{1}{2}\, \rho^{\veps}_{0,\alpha,T}(f)\, e^{\veps|x|} \int_{0}^t \,  |\frac{1-\tau}{1-t}|^\veps|\frac{1+\tau}{1+t}|^{-\veps} \frac{1}{|1-\tau|}  d\tau \nonumber \\
 & \leq & \frac{1}{2}\,\rho^{\veps}_{0,\alpha,T}(f)\, e^{\veps|x|} \veps^{-1} |1-t|^{-\veps}
\end{eqnarray}
By using the same arguments, when $-1<t<0,$ we obtain the same estimate to $G_+ f.$ Therefore
$$
|\partial^\alpha_x(G_+ f(t,x))| \ \leq \ \frac{1}{2}\,\rho^{\veps}_{0,\alpha,T}(f)\, e^{\veps|x|} \veps^{-1} |1+t|^{-\veps}
$$
Therefore, given $f\in C^\infty (\R: E^\veps (\R)), \alpha \in \Z_+$ and $K\subset\subset \R$,  we set
$$M = \veps^{-1}\, \rho^{\veps}_{0,\alpha,T}(f)\, \sup_{x \in K} e^{\veps|x|}.$$
Thus, it follows from \eqref{modvf3} and \eqref{modvf4} that
$$
|\partial^\alpha_x(L^{-1}_0 f(t,x))| \ \leq  M\, |1\pm t|^{-\veps}, \ \quad 0 < |1\pm t| \ll 1,\, x\in K
$$ \vskip-1.5em
\end{proof}
\

\begin{remark}
Since the general solution of $L_0 u=f$ in $[-1,1] \times \R$ is given by
\begin{equation}\label{example-sharp}
u = \varphi (x + \ln (1-t^2)) + L_0^{-1}f(t,x),
\end{equation}
with $\varphi$ being a function (or distribution) of one variable, we observe that we have always singularity at $t=-1$ or $t=+1$.

In view of the separatrix phenomena, we have not compensate both singularities in the general case, while we can ``cancel" the singularity either at $t=-1$ or $t=+1$

If $f \equiv c \neq 0,$ we exhibit, apart from  $u = \frac{c}{2}\ln|\frac{1+t}{1-t}|$, two particular solutions:
\begin{equation}
u_{\pm}(t,x) = \mp \frac{c}{2}x \pm c\ln|1\mp t|
\end{equation}
\end{remark}

\section{Perturbation with nondegenerate p.d.o.} \

The aim of this section is to show that if we perturb $L_0$ with constants p.d.o., or more generally, a Fourier multiplier satisfying suitable nondegeneracy conditions, we can obtain $L^\infty_{loc}$ estimates in $t$ for the $(L_0 + b)^{-1}$ without the smallness requirement on $b$.

More precisely, we consider 
  \begin{equation}
L_bu  =  (1-t^2)\partial_t u  - 2t \partial_xu + b(D)u=  f(t,x)
\label{vfm1}
\end{equation}
where
 \begin{equation}
\textrm{$b(\xi) \in C(\R)$ is realvalued and bounded away from zero for $\xi\in \R$.}
\label{mult1}
\end{equation}
Clearly \eqref{mult1} implies that one can find $0<\delta_0 <\delta_1$ such that
\begin{equation}
\textrm{either $\delta_0 \leq b(\xi)\leq \delta_1$ \ or \ $-\delta_1 \leq b(\xi) \leq -\delta_0$, for $\xi\in \R$.}
\label{mult2}
\end{equation}

Set $\hat{u}(t,\xi) = \mathcal{F}_{x\to \xi}u(t,\cdot)$ to be the partial Fourier transform in $x$, i.e.,

\begin{equation}
\hat{w}(\xi)= \int_{\R} e^{-ix\xi} w(x) dx.
\label{mult3}
\end{equation}

  Setting (formally)
   \begin{equation}
\hat{u} (t,\xi)   = \exp( -\frac{b(\xi)}{2}\ln \frac{|1+t|}{|1-t|}  ) =  \left( \frac{|1-t|}{|1+t|}\right)^{b(\xi)/2} \hat{v}(t,\xi)
\label{vf2b}
\end{equation}
we obtain that
 \begin{equation}
\hat{L_0}\hat{v}(t,\xi)    = \left( \frac{|1-t|}{|1+t|}\right)^{-b(\xi)/2} \hat{f}(t,\xi).
\label{mult4}
\end{equation}

In view of the nondegeneracy condition \eqref{mult2}  can write a right inverse of $L_b$ which is $L^\infty_{loc} $ in $t$ (a better regularity than $L^1_{loc}$ for $L^{-1}$). Indeed, set

 \begin{eqnarray}
\hat{L_b}^{-1} \hat{f} & = & \left( \frac{|1-t|}{|1+t|}\right)^{b(\xi)/2} \!\! \int_{-\sign (b)}^t \frac{\hat{f}(s, \xi) e^{i \ln \frac{|1-t^2|}{|1-\tau^2|}}}{(1-\tau) |1-\tau|^{b(\xi)/2} (1+s) |1+s|^{-b(\xi)/2}}\, ds \ \qquad
\label{vf4b}
\end{eqnarray} \\

\begin{proposition}\label{pert1const}
The operator $L_b^{-1}$ acts continuously as $L_0^{-1}$ in the spaces with subexponential decay. Furthermore, it acts continuously
 \begin{eqnarray}
L_b^{-1} & : & C(\R: H^s(\R) ) \longmapsto L^\infty_{loc} (\R: H^s(\R))
\label{Sobpert1}
\end{eqnarray}
and for every $K>0$, $s>0$, one can find $C=C_K>0$ such that
\begin{eqnarray}
\|L_b^{-1}f\|_{L^\infty ([-K,K]:H^s(\R))} & \leq & \frac{C}{\delta_0}\| f\|_{C([-K,K]: H^s(\R))},
\label{Sobpert2}
\end{eqnarray}
for all $f\in C(\R: H^s(\R) )$ and $\delta_0 >0$.
\end{proposition}

\begin{proof}

We have the crucial step is based on the estimates near $t=\pm 1$:
 \begin{eqnarray}
\| \hat{L_b^{-1}} \hat{f}(t,\cdot)\|_{L^2}
& \leq &  C_0 \|f\|_{C([-K,K]:L^2(\R))} \sup_{\xi\in \R} \left( |1+\sign (b) t|\right)^{\pm b(\xi)/2} \nonumber \\
&  &\times \ |\int_{-\sign (b)}^t \frac{1}{|1- s|^{1+ \sign (b)/2} |1+s|^{1- \sign (b)/2}} \, ds| \nonumber \\
& \leq & \frac{C}{|b|} \| f\|_{C([-K,K]:L^2(\R))}
\label{ppert1}
\end{eqnarray}
where 
$$C_0 = \sup_{\xi \in \R} \left( (\frac{|1-t|}{|1+t|})^{b(\xi)/2} \!\! \int_{-\sign (b)}^t \frac{1}{(1-\tau) |1-\tau|^{b(\xi)/2} (1+s) |1+s|^{-b(\xi)/2}}\, ds\right) \leq \frac{2}{\delta_0} 
$$
\end{proof}

\section{Final Remarks} \

First we observe that our results remain valid for vector fields of the type
$$
L = p(t)\partial_t + q(t,x) \partial_x
$$
provided $q(t,x)$ is bounded for $x$, when $x \to \infty$. The approach follows the same ideas, but the arguments of the proofs become more involved in view of the use of theorems on global behaviour of solutions of o.d.e. If $q$ is not bounded, for $x \to \infty$, we have more restrictive conditions on the growth of the RHS $f$. For example, if  $q(t,x)$ grows linearly in $x$ (like SG first order hyperbolic pseudodifferential operators (cf. \cite{Co1}), we have to require that the RHS $f(t,x)$ grows less than every $|x|^\gamma,$ for every $\gamma >0.$ Next, we point out that if the RHS $f$ decays to zero for $x\to \infty$, the right inverses $L_j$.

Next, as to possible multidimensional generalizations of the vector fields studied in the present work, we are also able to propose similar results for some classes of vector fields having smooth symmetries.
E.g. consider the regular plane vector field $ L=(t^2-15)(t^2+15)\partial_x-(t^2-25)(t^2-9)t\partial_t\,.$ One can easily check that the rotations of $L$ around the $x$ axis in $\R^3$ with coordinates $(t,x,y)$ gives rise to a regular vector field $M$ having as separatrices the two cilinders $y^2 + t^2 = 9$ and $y^2 + t^2 = 25$. The cohomological equation $M u=v$ hence is not solvable for every smooth function $v\in C^\infty(\R^3)$ because of the theorem of Duistermaat and Hormander but our techniques can be used to find weak solutions.

Finally, we point out to a natural  problem related to the reduction of a perturbation $L+b(t,x,D)$ to $L$ by means of global conjugation formally $J (t)\circ(L+b)\circ J^{-1}(t) = L$, with $J$ being a global p.d.o. or Fourier integral operator in $x\in \R^n$ depending smoothly on $t\in R\setminus I_L$, with singularities near $t=t_j$, $S_j$ or $S_{j+1}$ being separatrix strips. The example in Section 4 suggests that one should aim on estimates of $J(t)$ in $L^1_{loc} (\R: B(\R^n))$, where $B(\R^n)$ stands for some weighted Sobolev type space (cf. \cite{CNR1}, \cite{Ru1}, \cite{RuSu1}
and the references therein for global estimates in $\R^n$ for Fourier integral operators).

\subsection*{Acknowledgment}
The authors are grateful to Adalberto Bergamasco and Michael Ruzansky for useful discussions on arguments related to this article.

% ------------------------------------------------------------------------
\end{document}